\documentclass{article}
\pdfoutput=1

\usepackage[narrow]{arxiv}

\usepackage[utf8]{inputenc} % allow utf-8 input
\usepackage[T1]{fontenc}    % use 8-bit T1 fonts
\usepackage[hypertexnames=false]{hyperref}       % hyperlinks
\usepackage{url}            % simple URL typesetting
\usepackage{booktabs}       % professional-quality tables
\usepackage{amsfonts}       % blackboard math symbols
\usepackage{nicefrac}       % compact symbols for 1/2, etc.
\usepackage{microtype}      % microtypography

\usepackage{lipsum}         % Can be removed after putting your text content
\usepackage{graphicx}
\usepackage{natbib}
\usepackage{doi}

\title{An Optimal Second-Order Algorithm for Finite-Sum Optimization under Average Smoothness and Convexity}
\renewcommand{\shorttitle}{Second Order Finite-Sum Optimization under Average Smoothness and Convexity}

\newif\ifuniqueAffiliation
\uniqueAffiliationtrue

\ifuniqueAffiliation % Standard variant of author block
\author{
  Artem Tsedenov\\
  Yandex Research \\
  \texttt{cedenov18092005@gmail.com}
  \And
  Dmitry Kovalev \\
  Yandex Research \\
  \texttt{dakovalev1@gmail.com}
}
\else
\usepackage{authblk}

\author[1]{%
  Dmitry Kovalev\thanks{\texttt{dakovalev1@gmail.com}}%
}
\affil[1]{Yandex Research}
\fi
\renewcommand{\headeright}{}
\usepackage{algorithm,algpseudocodex}
\usepackage{nicematrix}

\usepackage[textsize=tiny,color=orange!20,bordercolor=orange,]{todonotes}
\presetkeys{todonotes}{inline}{}

\usepackage{xcolor}
\hypersetup{
  colorlinks=true,
  linkcolor=red!75!black,
  citecolor=blue!50!black
}

\usepackage{mathtools,amsthm,amssymb}
\allowdisplaybreaks[4]

\usepackage{enumitem}
\setlist{topsep=0pt,partopsep=0pt,itemsep=0pt}

\makeatletter\AddToHook{env/algorithmic/before}{\def\@currentcounter{ALG@line}}\makeatother
\usepackage{cleveref}
\crefalias{ALG@line}{line}
\usepackage{annotate}
\usepackage{mythm}
\usepackage{mymath}

\newcommand{\Rd}{\R^d}
\newcommand{\ox}{\overline{x}}
\newcommand{\hx}{\hat{x}}

\newcommand{\heta}{\hat{\eta}}
\newcommand{\N}{\mathbb{N}}
\newcommand{\ind}[1]{\mathbb{I}_{\brf{#1}}}
\newcommand{\be}{\mathbf{e}}
\newcommand{\taylor}{\operatorname{Taylor}}

\begin{document}
\maketitle

\begin{abstract}
  In this paper, we consider a finite-sum optimization problem, where the objective function is convex, the $n$ component functions are twice continuously differentiable and their Hessians are mean-square Lipschitz. Our contribution is twofold. First, we develop a stochastic algorithm, which requires $\tilde{\cO}(n + n^{6/7}/\epsilon^{2/7})$ second-order oracle calls in expectation to find an approximate solution to the problem with the expected accuracy $\epsilon$. Second, we establish a lower complexity bound of $\Omega(n + n^{6/7}/\epsilon^{2/7})$ for any stochastic second-order algorithm satisfying a stochastic second-order version of the first-order linear span assumption, which is widely adopted in the optimization literature. Our lower and upper complexity bounds match up to logarithmic factors and thus resolve an important open question of establishing the optimal complexity of this problem class.
\end{abstract}

\section{Introduction}

We consider the following finite-sum optimization problem:
\begin{equation}\label{eq:main}
  \min_{x \in \Rd} \brs*{ f(x) = \tmfrac{1}{n}\msum_{i=1}^n f_i(x)},
\end{equation}
where the objective function $f(x)\colon \Rd \to \R$ is convex, and the functions $f_i(x)\colon \Rd \to \R$ are twice continuously differentiable. We assume that this problem has a solution $x^* \in \Rd$. We also assume that the Hessians $\nabla^2 f_i(x)$ are mean-square Lipschitz, that is, the following condition holds:
\begin{aequation}<ass:L>
  \tmfrac{1}{n}\msum_{i=1}^n \sqn{\nabla^2 f_i(x) - \nabla^2 f_i(x')} \leq L^2 \sqn{x-x'}
  \quad
  \text{for all}\; x,x' \in \Rd.
\end{aequation}
The main goal of this paper is to develop an {\em optimal stochastic second-order algorithm} for finding an approximate solution to problem~\eqref{eq:main}.

By {\em second-order}, we mean that the algorithm can compute the gradients and the Hessians of the functions $f(x)$ and $f_i(x)$ at any given point $x \in \Rd$. Such computations are typically referred to as {\em second-order oracle calls}. By {\em stochastic}, we mean that the algorithm has access to individual functions $f_i(x)$ rather than the full average $f(x)$, where the index $i$ is selected according to some possibly randomized rule. Consequently, the computation of the full gradient $\nabla f(x)$ and Hessian $\nabla^2 f(x)$ requires $n$ calls to the oracles associated with the functions $f_i(x)$, and we aim to use as few such calls as possible. Finally, by {\em optimal}, we mean that we seek to establish a theoretical upper bound on the number of second-order oracle calls required by the algorithm to find an $\epsilon$-accurate solution\footnote{A point $\hx \in \Rd$ is called an $\epsilon$-accurate solution to problem~\eqref{eq:main} if it satisfies $\E{f(\hx) - f(x^*)} \leq \epsilon$.} to problem~\eqref{eq:main}, as well as a matching lower bound.

\subsection{Related Work}\label{sec:related}

\textbf{Deterministic first-order methods.} Standard gradient descent requires $\cO(n/\epsilon)$ gradient computations to find an $\epsilon$-accurate solution to a convex problem with a Lipschitz continuous gradient \citep{nesterov2013introductory}. This gradient complexity can be improved using Nesterov acceleration \citep{nesterov2013introductory,nesterov1983method}: accelerated gradient descent achieves a complexity of $\cO(n/\sqrt{\epsilon})$. The factor of $n$ appears in both cases because computing the full gradient $\nabla f(x)$ requires evaluating all $n$ individual gradients $\nabla f_i(x)$. The latter complexity bound cannot be improved further, which is implied by
the existing lower bounds \citep{nesterov2013introductory}.

\textbf{Deterministic second-order methods.} Standard Newton's method can diverge on convex problems with Lipschitz Hessians \citep{vandenberghe2004convex}. Global convergence can be achieved by Newton's method with cubic regularization \citep{nesterov2006cubic}: this algorithm achieves a second-order oracle complexity of $\cO(n/\sqrt{\epsilon})$. Cubic regularization can also be combined with Nesterov acceleration \citep{nesterov2008accelerating}, which leads to the accelerated cubic-regularized Newton's method with a complexity of $\cO(n / \epsilon^{1/3})$. A further improvement is obtained using the accelerated Newton's proximal extragradient method, which is based on the A-HPE framework and achieves a complexity of $\tilde{\cO}(n / \epsilon^{2/7})$ \citep{monteiro2013accelerated}. The extra logarithmic factors in this bound were later removed by \citet{kovalev2022first,carmon2022optimal}, yielding a complexity of $\cO(n/\epsilon^{2/7})$, which is optimal, as shown by the existing lower bounds \citep{nesterov2021implementable,arjevani2019oracle}.

\textbf{Stochastic first-order methods.} Standard stochastic gradient descent can be applied to smooth convex finite-sum problems. However, it achieves a poor gradient complexity of $\cO((\sigma/\epsilon)^2)$ and requires that the variance of the stochastic gradient is bounded by a constant $\sigma^2$ \citep{lan2012optimal}. Convergence without the bounded-variance assumption can be achieved by stochastic gradient methods with variance reduction \citep{johnson2013accelerating,schmidt2017minimizing}. Such methods can achieve a gradient complexity of $\tilde{\cO}(n + 1/\epsilon)$ \citep{allen2016improved}, offering a significant improvement over full-batch gradient descent. A further improvement is achieved by accelerated stochastic  variance-reduced gradient methods \citep{lan2018optimal,allen2018katyusha,lan2019unified}, which have a complexity of $\tilde{\cO}(n + \sqrt{n/\epsilon})$, matching the existing lower bounds up to logarithmic factors \citep{woodworth2016tight}.

The latter result has an important caveat: it requires each individual function $f_i(x)$ to be convex and have a Lipschitz continuous gradient. These requirements can be relaxed by assuming that only the full average function $f(x)$ is convex and the individual gradients $\nabla f_i(x)$ are mean-square Lipschitz. In this case, the optimal complexity is $\cO(n + n^{3/4} / \sqrt{\epsilon})$. It was achieved (up to logarithmic factors) by the algorithms of \citet{allen2018katyushax,shalev2016sdca}, and a matching lower bound was established by \citet{zhou2019lower}. It is also important to highlight that the latter setting is the {\em direct first-order analog of the second-order setting of this paper}. Indeed, we only require a convex average function $f(x)$ and mean-square Lipschitz individual Hessians $\nabla^2 f_i(x)$ (\Cref{ass:L}).

\textbf{Stochastic second-order methods.} Several works analyze the convergence of second-order methods in the fully stochastic or inexact settings \citep{agafonov2024inexact,antonakopoulos2022extra,agafonov2023advancing}. Unfortunately, their stochastic guarantees require additional assumptions on the gradient and/or Hessian noise, do not fully exploit the finite-sum structure, and typically do not improve upon the $\cO((\sigma/\epsilon)^2)$ complexity of first-order SGD in the high-accuracy setting. There is also an array of works analyzing stochastic variance-reduced second-order methods in the {\em non-convex} setting \citep{zhou2019stochastic,wang2019stochastic,shen2019stochastic,zhou2020stochastic,chayti2023unified,wu2026matching}. This neighboring research direction is relatively well-studied: weak assumptions make it impossible to guarantee fast convergence to a global minimum, and state-of-the-art results can be achieved with relatively simple algorithms without requiring unbiased stochastic updates and without advanced mechanisms such as the extragradient method or Nesterov acceleration. On the other hand, existing results on stochastic variance-reduced second-order methods for convex functions are highly limited. \citet{chayti2023unified} guarantee a complexity of $\tilde{\cO}(n + n^{4/5}/\sqrt{\epsilon})$ for gradient-dominated functions, although they do not provide corresponding guarantees for convex functions.\footnote{The gradient domination property with exponent $1$ holds for convex functions on bounded sets. However, \citet{chayti2023unified} do not guarantee that the iterates are bounded.} There are also several local convergence results for strongly convex functions \citep{kovalev2019stochastic,derezinski2022stochastic}.

\subsection{Main Contributions}

In \Cref{sec:alg}, we develop an accelerated stochastic variance-reduced second-order method (\Cref{alg}) and prove that it requires $\cO(n\log n + n^{6/7}/\epsilon^{2/7})$ second-order oracle calls in expectation to find an $\epsilon$-accurate solution to problem~\eqref{eq:main} when the objective function $f(x)$ is convex and \Cref{ass:L} holds. In \Cref{sec:lb}, we prove a lower bound of $\Omega(n + n^{6/7}/\epsilon^{2/7})$ on the expected number of second-order oracle calls required by any algorithm satisfying \Cref{ass:alg} to reach the same accuracy. These bounds coincide in the $\epsilon$-dependent term and match overall up to a logarithmic factor. Thus, \Cref{alg} is an optimal algorithm in the high-accuracy regime and optimal up to logarithmic factors overall.

\textbf{Notations.} In this paper, we will use the following notation:
$\norm{}_p$ denotes the standard vector $p$-norm for $p \in [1,\infty]$;
for a vector $x \in \Rd$, $\norm{x} = \norm{x}_2 = \sqrt{\<x,x>}$ denotes the standard Euclidean norm induced by the standard inner product $\<,>$;
for a matrix $\mA \in \R^{p\times q}$, $\norm{\mA} = \sup_{\norm{z}\leq1} \norm{\mA z}$ denotes the spectral norm;
$\N$ and $\N_0$ denote the sets of positive and nonnegative integers, respectively;
$\be_{d,j}\in \Rd$ denotes the $j$-th standard unit basis vector in $\Rd$ for $j \leq d$, and $\be_{d,j} = 0$ for $j > d$;
by $\taylor_{h}(x;w) = h(w) + \<\nabla h(w),x-w> + \frac{1}{2}\smash{\sqn{x-w}_{\nabla^2 h(w)}}$, we denote the second-order Taylor approximation of a twice continuously differentiable function $h(x)\colon \Rd \to \R$ at the point $w \in \Rd$.

\section{Optimal Stochastic Second-Order Algorithm}\label{sec:alg}

\subsection{The Algorithm and its Construction}

In this section, we describe an optimal stochastic second-order algorithm for solving problem~\eqref{eq:main}. Below, we briefly outline the main tools and components that we use to construct the algorithm.

\textbf{Extragradient method.} The extragradient method \citep{korpelevich1976extragradient} is the starting point for our algorithm. This is a natural choice because {\bf (i)} the existing optimal deterministic second-order methods are based on it \citep{monteiro2013accelerated,kovalev2022first,carmon2022optimal}, and {\bf(ii)} it has a stochastic extension, extragradient with variance reduction \citep{alacaoglu2022stochastic}, which has convergence guarantees for finite-sum convex problems with mean-square Lipschitz gradients. The mean-square Lipschitz gradient assumption is the first-order analog of \Cref{ass:L}, as discussed in \Cref{sec:related}.

\textbf{Newton's method with cubic regularization.} The standard extragradient method can be combined with the cubic-regularized Newton's method \citep{kovalev2022first}. The iterations of the first-order and the second-order variants of the extragradient method are given below in \cref{eq:eg1,eq:eg2}, respectively:
\begin{align}
  x_{k+1} &= x_k - \eta_k \nabla f(x_{k+1/2}),
  &
  x_{k+1/2} &= x_k - \eta_k \nabla f(x_k),\label{eq:eg1}
  \tag{EG.1}
  \\
  x_{k+1} &= x_k - \eta_k \nabla f(x_{k+1/2}),
  &
  x_{k+1/2} &= x_k - \eta_k \nabla \taylor_{f}(x_{k+1/2};x_k),%\brs{\nabla f(x_k) + \nabla f^2(x_k)(x_{k+1/2} - x_k)},
  \label{eq:eg2}
  \tag{EG.2}
\end{align}
where $\eta_k > 0$ is constant in \cref{eq:eg1}, and $\eta_k =  [L\norm{x_{k+1/2} - x_k}]^{-1}$ in \cref{eq:eg2}. The update rule for $x_{k+1/2}$ in \cref{eq:eg2} is well-defined and can be expressed as a single step of Newton's method with cubic regularization starting from $x_k$ or, equivalently, as a solution to the following auxiliary optimization problem:
\begin{equation}
  x_{k+1/2} = \targmin_{x \in \Rd}\brs{\taylor_f(x;x_k) + \tfrac{1}{3}L\norm{x-x_k}^3}.
\end{equation}

\textbf{Variance reduction.}  We now turn to a stochastic variance-reduced extension of the deterministic extragradient method \citep{alacaoglu2022stochastic}. This extension uses the following updates:
\begin{align}
  x_{k+1/2} &= x_k + \pi_k(w_k - x_{k+1/2}) -  \eta_k \nabla f(w_k),\label{eq:veg1a}  \tag{vEG.1a}
  \\
  x_{k+1} &= x_{k} + \pi_k(w_k - x_{k+1/2}) - \eta_k \brs{\nabla f_{i_k}(x_{k+1/2}) - \nabla f_{i_k}(w_k) + \nabla f(w_k)},\label{eq:veg1b}  \tag{vEG.1b}
\end{align}
where $\pi_k\in(0,1)$ and $\eta_k > 0$ are constant, $w_k \in \Rd$ is a reference point, and $i_k \sim \text{u.a.r}(\brf{1,\dots,n})$. The reference point $w_k$ is updated with probability $1/n$ at each iteration; thus, the algorithm requires $\cO(1)$ gradient evaluations per iteration on average. Motivated by the parallel between the deterministic first- and second-order extragradient methods above, we propose the following stochastic variance-reduced second-order extragradient method:
\begin{align}
  x_{k+1/2} &= x_k + \pi_k(w_k - x_{k+1/2}) -  \eta_k \nabla \taylor_f(x_{k+1/2}; w_k)
  \label{eq:veg2a}\tag{vEG.2a}
  \\
  x_{k+1} &= x_{k} + \pi_k(w_k - x_{k+1/2}) - \eta_k \brs{\nabla f_{i_k}(x_{k+1/2}) + \nabla \taylor_{f - f_{i_k}}(x_{k+1/2}; w_k)},\label{eq:veg2b}\tag{vEG.2b}
\end{align}
where $\eta_k = \sqrt{\pi_k} / [L\norm{x_{k+1/2} - w_k}]$. Unfortunately, the update rule in \cref{eq:veg2a} no longer corresponds to the solution of an auxiliary optimization problem with cubic regularization. However, we will show below that it is sufficient to use a clipped version of this stepsize, which can be computed using a simple line-search procedure. The update rules in \cref{eq:veg1b,eq:veg2b} each use an unbiased stochastic estimator of the gradient $\nabla f(x_{k+1/2})$. The update rule in \cref{eq:veg2b} is also related to the SVRG2 algorithm of \citet{gower2018tracking}: they use a similar update rule of the form
\begin{equation}
  x_{k+1} = x_{k} - \eta_k \brs{\nabla f_{i_k}(x_k) + \nabla \taylor_{f - f_{i_k}}(x_k; w_k)}
\end{equation}
with a constant $\eta_k$. However, \citet{gower2018tracking} establish a global convergence rate, which does not improve upon the rate of first-order SVRG, and only for generalized linear models under restrictive assumptions. Consequently, their results establish neither acceleration nor optimal convergence rates in the stochastic finite-sum second-order optimization setting.

\textbf{Nesterov acceleration.} While the first-order extragradient method with variance reduction enjoys optimal convergence rates for a more general class of finite-sum monotone inclusion problems, achieving optimal guarantees for minimization problems requires incorporating Nesterov acceleration. One of the simplest ways to incorporate acceleration is via the framework of \citet{borodich2025linear}, which was also used by \citet{xie2026tale,kovalev2026sgd,kovalev2026optimal} for other problems. The idea is to apply variance-reduced extragradient to the modified functions $f_k(x),f_{k,i}(x)\colon \Rd \to \R$, which are defined as follows:
\begin{equation}\label{eq:fk}
  f_k(x) = \tmfrac{1}{n} \msum_{i=1}^n f_{k,i}(x)
  ,\qquad
  f_{k,i}(x) = \alpha_k^{-3} f_i(\alpha_k x + (1-\alpha_k) \ox_k),
\end{equation}
where $\ox_k \in \Rd$ and $\alpha_k \in (0,1]$. The Hessians $\nabla^2 f_{k,i}(x)$ remain mean-square Lipschitz with the same constant $L$. Thus, we can apply the convergence analysis of the non-accelerated method and use it to determine the update rules for $\ox_k$ and $\alpha_k$. Updating $\alpha_k$ and $\ox_k$ triggers the computation of the full gradient $\nabla f\brr{}$ and full Hessian $\nabla^2 f\brr{}$. We therefore aim to update $\ox_k$ and $\alpha_k$ only $\cO(1)$ times per $n$ iterations on average.

Based on the above discussion, we now complete the construction of the algorithm. For compactness, we define the second-order Taylor approximations $\cT_k(\cdot;w),\cT_{k,i}(\cdot;w)\colon \Rd \to \R$  of the functions $f_k(x),f_{k,i}(x)$ as follows:
\begin{equation}\label{eq:taylor}
  \cT_{k}(x;w) = \tmfrac{1}{n}\msum_{i=1}^n\cT_{k,i}(x;w)
  ,\quad
  \cT_{k,i}(x;w) = \taylor_{f_{k,i}}(x;w).
\end{equation}
We also define auxiliary functions $\hx_{k,t}(\eta) \colon \R_+ \to \Rd$ and $\psi_{k,t}(\eta)\colon \R_+\to \R$ as follows:
\begin{equation}\label{eq:hx}
  \begin{aligned}
    \hx_{k,t}(\eta) &= \argmin_{x\in \Rd} \brs*{\tfrac{1}{2}\sqn{x-x_{k,t}} + \tfrac{1}{2}\pi_k\sqn{x-w_{k,t}} + \eta \cT_k(x;w_{k,t})},\\
    \psi_{k,t}(\eta) &= \eta L\norm{\hx_{k,t}(\eta) - w_{k,t}},
  \end{aligned}
\end{equation}
where $x_{k,t},w_{k,t} \in \Rd$ and $\pi_k > 0$. Finally, we present \Cref{alg}, which has two nested loops. The main purpose of the outer $k$-loop is to update $\ox_k$ and $\alpha_k$ to achieve accelerated convergence. The inner $t$-loop performs the iterations of the stochastic variance-reduced second-order extragradient method proposed above. It also incorporates the auxiliary line-search procedure given in \Cref{alg:ls}. The goal of this procedure is to find an approximate solution $\eta_{k,t}$ to the equation $\psi_{k,t}(\eta_{k,t}) = \sqrt{\pi_k}$, which corresponds to the desired stepsize $\eta_{k,t} \approx  \sqrt{\pi_k} / [L\norm{x_{k,t+1/2} - w_{k,t}}]$. This procedure also clips the stepsize $\eta_{k,t}$ on the final iteration of the $k$-th inner loop so that the cumulative sum of its stepsizes $\eta_{k,t}$ equals $\eta_k$ exactly, where $\brf{\eta_k}_{k \in \N_0} \subset \R_{++}$ is a predefined schedule. This is also necessary to guarantee accelerated convergence.

\begin{algorithm}[t]
  \caption{}
  \label{alg}
  \begin{algorithmic}[1]
    \State {\bf input:} $x_0=\ox_0 = w_0 = 0 \in \Rd$
    \State {\bf parameters:} $K \in \N$, $\brf{(\alpha_k,\pi_k,\eta_k)}_{k=0}^{K-1} \subset (0,1]^2 \times \R_{++}$
    \For{$k = 0,\dots,K-1$}
    \State $(x_{k,0},w_{k,0}) = (x_k,w_k)$
    \For{$t=0,1,2,\dots$}
    \State $\heta_{k,t} = \eta_k - \sum_{s=0}^{t-1}\eta_{k,s}$
    \If{$\heta_{k,t} = 0$}\label{line:if}
    \State $T_k = t$, \textbf{break}
    \EndIf
    \State $\eta_{k,t}= \Call{LineSearch}{k,t}$\label{line:eta_kt} \Comment{apply \Cref{alg:ls}}
    \State $x_{k,t+1/2} = \hx_{k,t}(\eta_{k,t})$\label{line:x_kt+/2} \Comment{$\hx_{k,t}(\eta)$ is defined in \cref{eq:hx}}
    \State sample $i_{k,t} \sim \text{u.a.r.}(\brf{1,\dots,n})$\label{line:i}
    \State sample $\beta_{k,t} \sim \text{Be}(\pi_k)$\label{line:beta_kt}
    \State $g_{k,t} = \nabla f_{k,i_{k,t}}(x_{k,t+1/2}) - \nabla \cT_{k,i_{k,t}}(x_{k,t+1/2};w_{k,t}) + \nabla \cT_k(x_{k,t+1/2};w_{k,t})$\label{line:g}
    \State $x_{k,t+1} = x_{k,t} + \pi_k(w_{k,t} - x_{k,t+1/2}) - \eta_{k,t} g_{k,t}$\label{line:x_kt+}
    \State $w_{k,t+1} = (1-\beta_{k,t})w_{k,t} + \beta_{k,t} x_{k,t+1/2}$ \label{line:w}
    \EndFor
    \State $(x_{k+1},w_{k+1}) = (x_{k,T_k},w_{k,T_k})$\label{line:xw_k+}
    \State $x_{k+1/2} = (1/\eta_k)\sum_{t=0}^{T_k-1}\eta_{k,t}x_{k,t+1/2}$\label{line:x_k+/2}
    \State sample $\beta_k \sim \text{Be}(\pi_k)$\label{line:beta_k}
    \State $\ox_{k+1} = (1-\beta_k)\ox_k + \beta_k(\alpha_k x_{k+1/2} + (1-\alpha_k)\ox_k)$\label{line:ox}
    \EndFor
    % \Statex
  \end{algorithmic}
\end{algorithm}

\begin{algorithm}[t]
  \caption{}
  \label{alg:ls}
  \begin{algorithmic}[1]
    \Procedure{LineSearch}{$k,t$}
    \If{$\psi_{k,t}^2(\heta_{k,t}) \leq \tfrac{3}{2}\pi_k$}\Comment{$\psi_{k,t}(\eta)$ is defined in \cref{eq:hx}}
    \State\Return $\heta_{k,t}$
    \EndIf
    \State $(\eta_{k,t,0}^-, \eta_{k,t,0}^+) = (0, \heta_{k,t})$
    \For{$l=0,1,2,\dots$}
    \State $\eta_{k,t,l} = \tfrac{1}{2}(\eta_{k,t,l}^- + \eta_{k,t,l}^+)$
    \If{$\psi_{k,t}^2(\eta_{k,t,l}) > \tfrac{3}{2}\pi_k$}
    \Comment{$\psi_{k,t}(\eta)$ is defined in \cref{eq:hx}}
    \State $(\eta_{k,t,l+1}^-, \eta_{k,t,l+1}^+) = (\eta_{k,t,l}^-, \eta_{k,t,l})$
    \ElsIf{$\psi_{k,t}^2(\eta_{k,t,l}) < \pi_k$}
    \Comment{$\psi_{k,t}(\eta)$ is defined in \cref{eq:hx}}
    \State $(\eta_{k,t,l+1}^-, \eta_{k,t,l+1}^+) = (\eta_{k,t,l}, \eta_{k,t,l}^+)$
    \Else
    \State\Return $\eta_{k,t,l}$
    \EndIf
    \EndFor
    \EndProcedure
  \end{algorithmic}
\end{algorithm}

\subsection{Convergence Analysis}

In this section, we analyze the convergence of \Cref{alg}. We start with \Cref{lem:eta_kt}, which establishes that \Cref{alg:ls} indeed finds an approximate solution to the equation $\psi_{k,t}(\eta_{k,t}) = \sqrt{\pi_k}$.\footnote{The final iteration $t = T_k - 1$ is excluded because the stepsize may be clipped.} The proof in \Cref{proof:lem:eta_kt} also shows that {\bf (i)} each execution of \Cref{alg:ls} terminates after finitely many iterations, {\bf (ii)} each iteration requires only basic linear algebra operations such as matrix-vector multiplications and matrix inversions, and {\bf (iii)} the line-search procedure requires no additional second-order oracle calls because it uses only the previously computed derivatives $\nabla f_k(w_{k,t})$ and $\nabla^2 f_k(w_{k,t})$.

\begin{lemma}<lem:eta_kt>
  $\eta_{k,t} \in \R_{++}$ and $x_{k,t+1/2} \in \Rd$ are well-defined on \Cref{line:eta_kt,line:x_kt+/2}, respectively, and satisfy $\pi_k \leq \psi_{k,t}^2(\eta_{k,t})$ for all $0 \leq t < T_k-1$ and $\psi_{k,t}^2(\eta_{k,t}) \leq \frac{3}{2}\pi_k$ for all $0 \leq t < T_k$.
\end{lemma}

An important aspect of \Cref{alg} is that the length $T_k$ of the $k$-th inner loop is a random variable. Most such algorithms in the literature sample the loop length from a distribution fixed in advance (e.g., \citet{lei2017less,kovalev2020don,carmon2022recapp}). \Cref{alg} differs substantially from these methods: the inner-loop length $T_k$ is determined by the stopping condition on \Cref{line:if}, rather than independently sampled from a fixed distribution. Consequently, the inner loop may never terminate, and in this case we set $T_k = +\infty$. However, we show below that $T_k$ is almost surely finite for every $k$ and that the average inner-loop length is $\cO(1)$ in expectation.

For clarity, let $\cF_{k,t}$ denote the sigma-algebra generated by all random variables sampled before the $t$-th inner iteration of the $k$-th inner loop, including $(\beta_j,\beta_{j,s},i_{j,s})$ for $j < k$, and $(\beta_{j,s},i_{j,s})$ for $j = k$ and $s < t$. We also use the notation $\E[k,t]{\cdot} = \E{\cdot \given \cF_{k,t}}$ and define the distance function $\Phi_{k,t}$ as follows:
\begin{equation}\label{eq:Phi}
  \Phi_{k,t} = \tfrac{1}{2}\sqn{x_{k,t} - x^*} + \tfrac{1}{2}\sqn{w_{k,t} - x^*}.
\end{equation}
Next, we state \Cref{lem:one_step}, which provides a single-iteration analysis of the stochastic variance-reduced second-order extragradient method used in the inner loop.
\begin{lemma}<lem:one_step>
  Suppose that \Cref{alg} has reached the $t$-th inner iteration of the $k$-th inner loop and that $t < T_k$. Then, the following inequality holds:
  \begin{equation}
    \eta_{k,t}\brs*{f_k(x_{k,t+1/2}) - f_k(x^*)} + \tfrac{1}{4}\pi_k\sqn{x_{k,t+1/2} - w_{k,t}} \leq \Phi_{k,t} - \E[k,t]{\Phi_{k,t+1}}.
  \end{equation}
\end{lemma}

Based on \Cref{lem:one_step}, \Cref{lem:inner_trunc} provides an analysis of finitely many inner-loop iterations. Next, we establish that the inner-loop length $T_k$ is almost surely finite and provide the convergence analysis for an entire inner loop in \Cref{lem:inner}. The proof of this lemma uses Fatou's lemma to pass to the limit in the expectation bound. To apply Fatou's lemma, we need the lower bound established in \Cref{lem:minorant}.

\begin{lemma}<lem:inner_trunc>
  Suppose that \Cref{alg} almost surely reaches the $k$-th iteration. Suppose that $\E{\Phi_{k,0}}\in \R$ and $\E{f(\ox_k)} \in \R$. Then, for all $T \in \N$, $\E{\Phi_{k,\min\brf{T+1,T_k}}} \in \R$ and the following inequality holds:
  \begin{equation}
    \begin{aligned}
      \mi{14}
      \E*{\msum_{t=0}^{\min\brf{T,T_k-1}} \brr*{\eta_{k,t}\brs*{f_k(x_{k,t+1/2}) - f_k(x^*)} + \tfrac{1}{4}\pi_k\sqn{x_{k,t+1/2} - w_{k,t}}}}
      \\&\leq
      \E{\Phi_{k,0}} - \E{\Phi_{k,\min\brf{T+1,T_k}}}
      \leq
      \E{\Phi_{k,0}}.
    \end{aligned}
  \end{equation}
\end{lemma}

\begin{lemma}<lem:minorant>
  Suppose that \Cref{alg} reaches the $k$-th iteration. Then, the following inequality holds for all $T \in \N_0$:
  \begin{align*}
    \msum_{t=0}^{\min\brf{T,T_k-1}} \eta_{k,t}\brs*{f_k(x_{k,t+1/2}) - f_k(x^*)} \geq \eta_k(1-\alpha_k)\alpha_k^{-3}\brs{f(x^*) - f(\ox_k)}.
  \end{align*}
\end{lemma}

\begin{lemma}<lem:inner>
  Suppose that \Cref{alg} almost surely reaches the $k$-th iteration. Suppose that $\E{f(\ox_k)}\in \R$ and $\E*{\Phi_{k,0}} \in \R$.
  Then, $T_k \in \N$ almost surely, and the following inequality holds:
  \begin{equation}
    \E*{\Phi_{k,T_k}} + \eta_k\E*{f_k(x_{k+1/2}) - f_k(x^*)} + \E*{\msum_{t=0}^{T_k-1}\tfrac{1}{4}\pi_k\sqn{x_{k,t+1/2} - w_{k,t}}} \leq \E{\Phi_{k,0}}.
  \end{equation}
\end{lemma}

Based on the above results, we state the key \Cref{lem:induction}. Combining this lemma with standard inductive arguments, we establish the main convergence result for \Cref{alg} in \Cref{thm:alg}. \Cref{thm:alg} also specifies the parameter choices for \Cref{alg}. Note that the parameter $\alpha_k$ changes once per $\Omega(n)$ outer iterations and that the expected inner-loop length $\E{T_k}$ is $\cO(1)$ when averaged over the outer iterations, as discussed above. We also provide an explicit bound on the second-order oracle complexity in \Cref{cor:alg}.

\begin{lemma}<lem:induction>
  Under the conditions of \Cref{lem:inner}, suppose that the following inequality holds:
  \begin{equation}
    \E*{\Phi_{k,0}} + \eta_k(1-\alpha_k\pi_k)\pi_k^{-1}\alpha_k^{-3}\E*{f(\ox_k) - f(x^*)} + \E*{\cW_k} \leq \sqn{x^*},
  \end{equation}
  where $\cW_k = \sum_{l=0}^{k-1}\sum_{t=0}^{T_l-1}\tfrac{1}{4}\pi_l\sqn{x_{l,t+1/2} - w_{l,t}}$.
  Then, the following inequality holds:
  \begin{equation}
    \E*{\Phi_{k+1,0} + \eta_k\pi_k^{-1}\alpha_k^{-3}\brs{f(\ox_{k+1}) - f(x^*)} + \cW_{k+1}} \leq \sqn{x^*}.
  \end{equation}
\end{lemma}

\begin{theorem}<thm:alg>
  Let $L > 0$ and $x^* \neq 0$.
  Let the parameters of \Cref{alg} be defined as follows:
  \begin{equation}\label{eq:params}
    \alpha_k =
    \begin{cases}
      1 & k< n\\
      1/\floor{\gamma_k} & k \geq n
    \end{cases}
    ,\quad
    \pi_k =
    \begin{cases}
      1/(k+1) & k < n\\
      1/n & k \geq n
    \end{cases}
    ,\quad
    \eta_k =
    \begin{cases}
      \eta & k < n\\
      \eta\alpha_k^3\gamma_k^{3.5} & k \geq n
    \end{cases},
  \end{equation}
  where $\gamma_k = 1 + \frac{2}{7n}[k - n + 1]$ and $\eta = 1/[16\sqrt{n}L\norm{x^*}]$. Consequently, for an arbitrary precision $\epsilon > 0$, choosing the number of outer iterations $K$ as
  \begin{equation}\label{eq:K}
    K
    =
    \ceil*{\tfrac{7}{2}n^{5/7}\brs*{\tmfrac{\sqn{x^*}}{\eta \epsilon}}^{2/7}}
    =
    \ceil*{\tfrac{7}{2}n^{6/7}\brs*{\tmfrac{16L \norm{x^*}^3}{\epsilon}}^{2/7}}
  \end{equation}
  implies $\E{f(\ox_K)-f(x^*)}\leq \epsilon$ and $\E{\sum_{k=0}^{K-1} T_k} \leq 2K$.
\end{theorem}

\begin{corollary}<cor:alg>
  \Cref{alg} requires $\cO\brr*{n\log (n+1) + n^{6/7}\brs{L\norm{x^*}^3/\epsilon}^{2/7}}$ computations of $\nabla f_i \brr{}$ and $\nabla^2 f_i\brr{}$ in expectation to reach the precision $\E{f(\ox_K)-f(x^*)}\leq \epsilon$.
\end{corollary}

\section{Lower Complexity Bounds}\label{sec:lb}

\subsection{Lower Bounds}

In this section, we show that the upper bound on the second-order oracle complexity in \Cref{cor:alg} cannot be improved beyond small logarithmic factors. We do this by establishing a lower bound on the expected number of iterations required by a stochastic second-order algorithm to find an $\epsilon$-accurate solution to problem~\eqref{eq:main}. We first specify a class of algorithms to which the lower bound applies. Concretely, we assume that each algorithm starts at the origin, i.e., $x_0 = 0$, and that its iterate $x_k \in \Rd$ satisfies the following condition at every iteration $k \in \N$:
\begin{aequation}<ass:alg>
  x_k \in \msum_{j=0}^{k-1} \spanset*{\brf*{h \in \Rd : h = a\nabla f_{i_j}(x_j) + b \nabla^2 f_{i_j}(x_j) h \;\;\text{for}\;\; a,b \in \R}},
\end{aequation}
where the indices $\brf{i_j}_{j \in \N_0}$ are sampled independently and uniformly from $\brf{1,\dots,n}$. This algorithm class is inspired by Assumption~1 of \citet{nesterov2021implementable} and can be seen as its direct stochastic extension.

Our lower complexity bounds are based on a specific hard instance of problem~\eqref{eq:main}. To simplify its construction, for arbitrary $d,n\in \N$ and $i \in \brf{1,\dots,n}$, we define the functions $h_{d,n}(x), h_{d,n,i}(x) \colon \Rd \to \R$ as follows:
\begin{equation}
  h_{d,n}(x) = \msum_{i=1}^n h_{d,n,i}(x)
  ,\qquad
  h_{d,n,i}(x) = \tfrac{1}{3}\msum_{j \in \cJ_i} \abs{\<\be_{d,j+1} - \be_{d,j},x>}^3,
\end{equation}
where $\cJ_i = \brf{j \in \N : 1 \leq j \leq d , j\equiv i \bmod n}$. Next, we establish the basic properties of these functions in the following lemmas. \Cref{lem:lb_L} is used to ensure that the hard instance satisfies \Cref{ass:L}. \Cref{lem:lb_sol} is used to derive the actual lower bound.

\begin{lemma}<lem:lb_sol>
  Let $\rho > 0$. Then $\argmin_{x \in \Rd} \brs{h_{d,n}(x) - \rho\<\be_{d,1},x>} = \sqrt{\rho} (d,d-1,\dots,1) \in \Rd$.
\end{lemma}

\begin{lemma}<lem:lb_L>
  The Hessians $\nabla^2 h_{d,n,i}(x)$ are mean-square Lipschitz, i.e., the following inequality holds:
  \begin{equation}
    \tmfrac{1}{n}\msum_{i=1}^n \sqn{\nabla^2 h_{d,n,i}(x) - \nabla^2 h_{d,n,i}(x')}
    \leq
    \tmfrac{256}{n}\sqn{x-x'}
    \quad\text{for all}\; x,x' \in \Rd.
  \end{equation}
\end{lemma}

Based on the above lemmas, we present our lower complexity bound in \Cref{thm:lb}. The proof extends the standard inductive arguments used, for example, by \citet[Section~4]{nesterov2021implementable}. The lower bound in \Cref{thm:lb} matches the upper bound in \Cref{cor:alg} up to the logarithmic factor $\log n$. In particular, the $\epsilon$-dependent term in the complexity cannot be improved.

\begin{theorem}<thm:lb>
  Let $L,R,\epsilon > 0$ and $n \in \N$. Then there exists an instance of problem~\eqref{eq:main} satisfying \Cref{ass:L} with a unique solution $x^* \in \Rd$ satisfying $\norm{x^*} = R$, such that the expected number of iterations required by any stochastic second-order algorithm satisfying \Cref{ass:alg} to output an $\epsilon$-approximate solution to problem~\eqref{eq:main} is
  \begin{equation}
    \Omega\brr*{n + n^{6/7} \cdot \brs*{LR^3/\epsilon}^{2/7}}.
  \end{equation}
\end{theorem}

\subsection{Discussion}

\textbf{Algorithm class.} The lower complexity bound in \Cref{thm:lb} applies to a specific class of stochastic second-order algorithms satisfying \Cref{ass:alg}. There are several important aspects to consider.
\begin{enumerate}[label={\bf(\roman{*})}]
  \item \Cref{alg} may seem inconsistent with \Cref{ass:alg} because it computes full-batch derivatives $\nabla f\brr{}$ and $\nabla^2 f\brr{}$. However, these computations can be implemented within the model in \Cref{ass:alg} with only logarithmic overhead. Indeed, we can compute the full-batch derivatives in \Cref{alg} by repeatedly sampling indices $i \sim \text{u.a.r.}(\brf{1,\dots,n})$ until every index has been sampled. This requires $\cO(n\log(n+1))$ samples in expectation, compared with $\cO(n)$ evaluations in the deterministic implementation. Alternatively, we can avoid this logarithmic overhead by extending \Cref{ass:alg} and the arguments in \Cref{proof:thm:lb} to allow the algorithm to evaluate derivatives of either a randomly selected individual function $f_i(x)$ or the full average $f(x)$. Such an extension is straightforward; we omit it because it would add extra technical detail without providing further insight.
  \item The linear span condition in \Cref{ass:alg} can probably be removed completely by extending the arguments of \citet{arjevani2019oracle}. Such an extension would not improve the $\epsilon$-dependent term in \Cref{thm:lb}. It is also not clear whether this extension would provide insight into the optimality of the term $\cO(n\log(n+1))$ in \Cref{cor:alg}. We leave this investigation for future research.
\end{enumerate}

\textbf{Optimal and near-optimal methods.} Even in light of the above discussion, the upper and lower complexity bounds in \Cref{cor:alg} and \Cref{thm:lb} do not match in all accuracy regimes. However, these bounds match in the high-accuracy regime, where $\epsilon$ is sufficiently small. Second-order methods are typically used when high accuracy is required; otherwise, first-order methods can be used without incurring the higher per-iteration cost of second-order methods. Outside the high-accuracy regime, the mismatch is only logarithmic, which makes \Cref{alg} near-optimal. We leave the question of optimality across all accuracy regimes for future research.

\newpage

\bibliographystyle{apalike}
\bibliography{references}

\newpage

\appendix

\crefalias{section}{appendix}
\crefalias{subsection}{subappendix}
\crefalias{subsubsection}{subsubappendix}

\section{Proofs for Section~\ref{sec:alg}}

\proofsubsection{lem:eta_kt}
One can observe that $\hx_{k,t}(\eta)$ can be explicitly computed as follows:
\begin{equation}
  \hx_{k,t}(\eta) = [(1+\pi_k)\mI_d + \eta \nabla^2 f_k(w_{k,t})]^{-1}[x_{k,t} + (\pi_k \mI_d + \eta \nabla^2 f_k(w_{k,t})) w_{k,t} - \eta\nabla f_k(w_{k,t})].
\end{equation}
Consequently, one can observe that $\hx_{k,t}(\eta)$ is a continuous function. Moreover, by standard limiting arguments, \Cref{alg:ls} finishes after finitely many iterations, and the output satisfies the desired inequalities.\qed

\proofsubsection{lem:one_step}

We can upper-bound $\frac{1}{2}\E[k,t]{\sqn{x_{k,t+1} - x^*}}$ as follows:
\begin{align*}
  \mi{0.5}\tfrac{1}{2}\E[k,t]{\sqn{x_{k,t+1} - x^*}}
  \\&\leq
  \tfrac{1}{2}\sqn{x_{k,t} - x^*}
  + \tfrac{1}{2}\E[k,t]{\sqn{x_{k,t+1} - x_{k,t+1/2}}}
  +\<{\E[k,t]{x_{k,t+1} - x_{k,t}}},x_{k,t+1/2} - x^*>
  \\&\at{uses \cref{line:x_kt+}}=
  \tfrac{1}{2}\sqn{x_{k,t} - x^*}
  +\tfrac{1}{2}\pi_k\sqn{w_{k,t} - x^*}
  -\tfrac{1}{2}\pi_k\sqn{x_{k,t+1/2} - x^*}
  \\&
  -\eta_{k,t}\<{\E[k,t]{g_{k,t}}},x_{k,t+1/2} - x^*>
  +\tfrac{1}{2}\E[k,t]{\sqn{x_{k,t+1} - x_{k,t+1/2}}}
  -\tfrac{1}{2}\pi_k\sqn{x_{k,t+1/2} - w_{k,t}}
  \\&\at{uses \cref{line:g,line:i} and the definitions in \cref{eq:fk,eq:taylor}}=
  \tfrac{1}{2}\sqn{x_{k,t} - x^*}
  +\tfrac{1}{2}\pi_k\sqn{w_{k,t} - x^*}
  -\tfrac{1}{2}\pi_k\sqn{x_{k,t+1/2} - x^*}
  \\&
  -\eta_{k,t}\<\nabla f_k(x_{k,t+1/2}),x_{k,t+1/2} - x^*>
  +\tfrac{1}{2}\E[k,t]{\sqn{x_{k,t+1} - x_{k,t+1/2}}}
  -\tfrac{1}{2}\pi_k\sqn{x_{k,t+1/2} - w_{k,t}}
  \\&\at{uses the convexity of the function $f_k(x)$}\leq
  \tfrac{1}{2}\sqn{x_{k,t} - x^*}
  +\tfrac{1}{2}\pi_k\sqn{w_{k,t} - x^*}
  -\tfrac{1}{2}\pi_k\sqn{x_{k,t+1/2} - x^*}
  \\&
  -\eta_{k,t}\brs{f_k(x_{k,t+1/2}) - f_k(x^*)}
  +\tfrac{1}{2}\E[k,t]{\sqn{x_{k,t+1} - x_{k,t+1/2}}}
  -\tfrac{1}{2}\pi_k\sqn{x_{k,t+1/2} - w_{k,t}}
  \\&\at{uses \cref{line:x_kt+/2,line:x_kt+,line:i,line:g} and the first-order optimality conditions in \cref{eq:hx}}=
  \tfrac{1}{2}\sqn{x_{k,t} - x^*}
  +\tfrac{1}{2}\pi_k\sqn{w_{k,t} - x^*}
  -\tfrac{1}{2}\pi_k\sqn{x_{k,t+1/2} - x^*}
  -\eta_{k,t}\brs{f_k(x_{k,t+1/2}) - f_k(x^*)}
  \\&
  +\tmfrac{1}{2n}\eta_{k,t}^2\msum_{i=1}^n\sqn{\nabla f_{k,i}(x_{k,t+1/2}) - \nabla \cT_{k,i}(x_{k,t+1/2};w_{k,t})}
  -\tfrac{1}{2}\pi_k\sqn{x_{k,t+1/2} - w_{k,t}}
  \\&\at{uses \Cref{ass:L} and the standard integration arguments}\leq
  \tfrac{1}{2}\sqn{x_{k,t} - x^*}
  +\tfrac{1}{2}\pi_k\sqn{w_{k,t} - x^*}
  -\tfrac{1}{2}\pi_k\sqn{x_{k,t+1/2} - x^*}
  -\eta_{k,t}\brs{f_k(x_{k,t+1/2}) - f_k(x^*)}
  \\&
  +\tfrac{1}{8}\eta_{k,t}^2L^2\norm{x_{k,t+1/2} - w_{k,t}}^4
  -\tfrac{1}{2}\pi_k\sqn{x_{k,t+1/2} - w_{k,t}}
  \\&\at{uses \Cref{lem:eta_kt}}\leq
  \tfrac{1}{2}\sqn{x_{k,t} - x^*}
  +\tfrac{1}{2}\pi_k\sqn{w_{k,t} - x^*}
  -\tfrac{1}{2}\pi_k\sqn{x_{k,t+1/2} - x^*}
  -\eta_{k,t}\brs{f_k(x_{k,t+1/2}) - f_k(x^*)}
  \\&
  -\tfrac{1}{4}\pi_k\sqn{x_{k,t+1/2} - w_{k,t}},
\end{align*}
where \annotate.
Next, from \Cref{line:w,line:beta_kt} it follows that
\begin{align*}
  \E[k,t]{\sqn{w_{k,t+1} - x^*}} = (1-\pi_k)\sqn{w_{k,t} - x^*} + \pi_k\sqn{x_{k,t+1/2} - x^*}.
\end{align*}
Combining this with the previous inequality concludes the proof.\qed

\proofsubsection{lem:inner_trunc}

First, we make the following definitions:
\begin{align*}
  A_{k,t} &= \ind{t< T_k}\eta_{k,t}\brs*{f_k(x_{k,t+1/2}) - f_k(x^*)}
  &
  B_{k,t} &= \tfrac{1}{4}\ind{t< T_k}\pi_k\sqn{x_{k,t+1/2} - w_{k,t}}
  \\
  C_{k,t} &= \ind{t< T_k}\Phi_{k,t}
  &
  D_{k,t} &= \E[k,t]{\ind{t< T_k}\Phi_{k,t+1}}.
\end{align*}
Using \Cref{lem:one_step}, for all $t \in \N_0$, we obtain
\begin{align*}
  A_{k,t} + B_{k,t}
  &\leq
  C_{k,t} - \ind{t< T_k}\E[k,t]{\Phi_{k,t+1}}
  \at{uses the fact that the condition $t < T_k$ is fully determined by $\cF_{k,t}$}=
  C_{k,t} - D_{k,t},
\end{align*}
where \annotate. It is obvious, that $\E{B_{k,t}}\geq 0$ and $\E{C_{k,t}} \geq 0$. Moreover, we can lower-bound $\E{A_{k,t}}$ as follows:
\begin{align*}
  \E{A_{k,t}}
  &\at{uses the definition in \cref{eq:fk}}\geq
  \E{\ind{t< T_k}\eta_{k,t}\alpha_k^{-3}\brs*{f(x^*) - f(\alpha_k x^* + (1-\alpha_k)\ox_k)}}
  \\&\at{uses the convexity of $f(x)$}\geq
  \E{\ind{t< T_k}\eta_{k,t}(1-\alpha_k)\alpha_k^{-3}\brs*{f(x^*) - f(\ox_k)}}
  \\&\geq
  \eta_k(1-\alpha_k)\alpha_k^{-3}\E{f(x^*) - f(\ox_k)} > -\infty,
\end{align*}
where \annotate.
Hence, after assuming $\E{C_{k,t}} \in \R$, we can upper-bound $\E{C_{k,t+1}}$ as follows:
\begin{align*}
  \E{C_{k,t+1}}
  &=
  \E{\ind{t+1< T_k}\Phi_{k,t+1}}
  \leq
  \E{\ind{t< T_k}\Phi_{k,t+1}}
  \at{uses the tower rule}=
  \E{D_{k,t}}
  \\&\leq
  \E{C_{k,t} - A_{k,t} - B_{k,t}}
  \leq
  \E{C_{k,t} - A_{k,t}}
  \at{uses the fact that $\E{A_{k,t}} > -\infty$}<
  +\infty,
\end{align*}
where \annotate. Since, $\E{C_{k,0}} = \E{\Phi_{k,0}} \in \R$, we can prove by induction that $\E{C_{k,t}} \in \R$ and $\E{D_{k,t}} \in \R$ for all $t \in \N_0$. Hence, we get
\begin{align*}
  \E*{\msum_{t=0}^T \brr*{A_{k,t} + B_{k,t}}} \leq \msum_{t=0}^T\brr*{\E{C_{k,t}} - \E{D_{k,t}}}.
\end{align*}
Finally, we plug all the definitions and obtain the following:
\begin{align*}
  \mi{3}
  \E*{\msum_{t=0}^{\min\brf{T,T_k-1}} \brr*{\eta_{k,t}\brs*{f_k(x_{k,t+1/2}) - f_k(x^*)} + \tfrac{1}{4}\pi_k\sqn{x_{k,t+1/2} - w_{k,t}}}}
  \\&\leq
  \msum_{t=0}^T
  \brr*{\E{\ind{t< T_k}\Phi_{k,t}} - \E{\E[k,t]{\ind{t< T_k}\Phi_{k,t+1}}}}
  \\&\at{uses the tower rule}=
  \msum_{t=0}^T
  \brr*{\E{\ind{t< T_k}\Phi_{k,t}} - \E{\ind{t< T_k}\Phi_{k,t+1}}}
  \\&=
  \msum_{t=0}^T
  \brr*{\E{\ind{t< T_k}\Phi_{k,t}} -\E{\ind{t+1< T_k}\Phi_{k,t+1}} - \E{\ind{t+1 = T_k}\Phi_{k,t+1}}}
  \\&=
  \E{\Phi_{k,0}}
  -\E{\ind{T+1< T_k}\Phi_{k,T+1}}
  -\msum_{t=0}^T \E{\ind{t+1 = T_k}\Phi_{k,T_k}}
  \\&=
  \E{\Phi_{k,0}}
  -\E{\ind{T+1< T_k}\Phi_{k,T+1}}
  -\E{\ind{T+1 \geq T_k}\Phi_{k,T_k}}
  \\&=
  \E{\Phi_{k,0}}
  -\E{\ind{T+1< T_k}\Phi_{k,\min\brf{T+1,T_k}}}
  -\E{\ind{T+1 \geq T_k}\Phi_{k,\min\brf{T+1,T_k}}}
  \\&=
  \E{\Phi_{k,0}}
  -\E{\Phi_{k,\min\brf{T+1,T_k}}},
\end{align*}
where \annotate.\qed

\proofsubsection{lem:minorant}
Similar to the proof of \Cref{lem:inner_trunc} in \Cref{proof:lem:inner_trunc}, we can obtain
\begin{align*}
  f_k(x_{k,t+1/2}) - f_k(x^*) \geq (1-\alpha_k)\alpha_k^{-3}\brs*{f(x^*) - f(\ox_k)}.
\end{align*}
It remains to observe that $\sum_{t=0}^{\min\brf{T,T_k-1}}\eta_{k,t} \leq \eta_k$.\qed

\proofsubsection{lem:inner}

From \Cref{lem:eta_kt}, for all $0 \leq t < T_k-1$, it follows that
\begin{align*}
  \brs*{\tfrac{1}{4}\pi_k\sqn{x_{k,t+1/2} - w_{k,t}}}\brs*{2L\eta_{k,t}/\pi_k}^2 \geq 1.
\end{align*}
Hence, we can upper-bound $\min\brf{T_k-1,T+1}$ for all $T \in \N_0$ as follows:
\begin{align*}
  \min\brf{T_k-1,T+1}
  &\leq
  \msum_{t=0}^{\min\brf{T_k-2,T}} \brs*{2L\eta_{k,t}/\pi_k}^{\frac{2}{3}}\brs*{\tfrac{1}{4}\pi_k\sqn{x_{k,t+1/2} - w_{k,t}}}^{\frac{1}{3}}
  \\&\at{uses H\"older's inequality}\leq
  \brs*{\msum_{t=0}^{\min\brf{T_k-2,T}} 2L\eta_{k,t}/\pi_k}^{\frac{2}{3}}
  \brs*{\msum_{t=0}^{\min\brf{T_k-2,T}} \tfrac{1}{4}\pi_k\sqn{x_{k,t+1/2} - w_{k,t}}}^{\frac{1}{3}}
  \\&\leq
  \brs*{\brs*{2L\eta_k/\pi_k}^2\msum_{t=0}^{\min\brf{T_k-2,T}} \tfrac{1}{4}\pi_k\sqn{x_{k,t+1/2} - w_{k,t}}}^{\frac{1}{3}},
\end{align*}
where \annotate. This implies the following:
\begin{align*}
  \E{\min\brf{T_k-1,T+1}^3}
  &\leq
  \brs*{2L\eta_k/\pi_k}^2\E*{\msum_{t=0}^{\min\brf{T_k-2,T}} \tfrac{1}{4}\pi_k\sqn{x_{k,t+1/2} - w_{k,t}}}
  \\&\at{uses \Cref{lem:inner_trunc,lem:minorant}}\leq
  \brs*{2L\eta_k/\pi_k}^2\E*{\Phi_{k,0} + \eta_k(1-\alpha_k)\alpha_k^{-3}\brs{f(\ox_k) - f(x^*)}}
  \\&<
  +\infty,
\end{align*}
where \annotate. By the monotone convergence theorem, we get the following:
\begin{align*}
  \E{(T_k-1)^3} < +\infty,
\end{align*}
which implies $T_k \in \N$ almost surely.
Finally, we use the definition
\begin{align*}
  A_{k,t} = \eta_{k,t}\brs*{f_k(x_{k,t+1/2}) - f_k(x^*)} + \tfrac{1}{4}\pi_k\sqn{x_{k,t+1/2} - w_{k,t}}
\end{align*}
and obtain the following inequality:
\begin{align*}
  \mi{1}\E{\Phi_{k,0}}
  \at{uses \Cref{lem:inner_trunc}}\geq
  \liminf_{T\to\infty}
  \E*{\Phi_{k,\min\brf{T+1,T_k}} + \msum_{t=0}^{\min\brf{T,T_k-1}} A_{k,t}}
  \\&\at{uses Fatou's lemma and \Cref{lem:minorant}}\geq
  \E*{\Phi_{k,T_k} + \msum_{t=0}^{T_k-1} A_{k,t}}
  \\&\at{uses \Cref{lem:minorant}}=
  \E*{\Phi_{k,T_k}} + \E*{\msum_{t=0}^{T_k-1}\brr*{\eta_{k,t}\brs*{f_k(x_{k,t+1/2}) - f_k(x^*)}}} + \E*{\msum_{t=0}^{T_k-1}\tfrac{1}{4}\pi_k\sqn{x_{k,t+1/2} - w_{k,t}}}
  \\&\at{uses \Cref{line:x_k+/2} and the convexity of the function $f_k(x)$}\geq
  \E*{\Phi_{k,T_k}} + \eta_k\E*{f_k(x_{k+1/2}) - f_k(x^*)} + \E*{\msum_{t=0}^{T_k-1}\tfrac{1}{4}\pi_k\sqn{x_{k,t+1/2} - w_{k,t}}},
\end{align*}
where \annotate.\qed

\proofsubsection{lem:induction}
We can upper-bound $\alpha_k^{-3}\E{f(\ox_{k+1}) - f(x^*)}$ as follows:
\begin{align*}
  \mi{3}\alpha_k^{-3}\E{f(\ox_{k+1}) - f(x^*)}
  \\&\at{uses \Cref{line:ox}}=
  \alpha_k^{-3}(1-\pi_k)\E{f(\ox_k) - f(x^*)}
  +\alpha_k^{-3}\pi_k\E{f(\alpha_k x_{k+1/2} + (1-\alpha_k)\ox_k) - f(x^*)}
  \\&\at{uses the definition in \cref{eq:fk}}=
  \alpha_k^{-3}(1-\pi_k)\E{f(\ox_k) - f(x^*)}
  \\&
  +\pi_k\E{f_k(x_{k+1/2}) - f_k(x^*) + \alpha_k^{-3}f(\alpha_k x^* + (1-\alpha_k)\ox_k) - \alpha_k^{-3}f(x^*)}
  \\&\at{uses the convexity of the function $f(x)$}\leq
  \alpha_k^{-3}(1-\pi_k)\E{f(\ox_k) - f(x^*)}
  \\&
  +\pi_k\E{f_k(x_{k+1/2}) - f_k(x^*) + (1-\alpha_k)\alpha_k^{-3}\brs{f(\ox_k) - f(x^*)}}
  \\&\at{uses the assumption $\E{f(\ox_k)}\in\R$}=
  \alpha_k^{-3}(1-\alpha_k\pi_k)\E{f(\ox_k) - f(x^*)}
  +\pi_k\E{f_k(x_{k+1/2}) - f_k(x^*)},
\end{align*}
where \annotate.
Furthermore, we can lower-bound $\E{\Phi_{k,0}}$ as follows:
\begin{align*}
  \E{\Phi_{k,0}}
  &\at{uses \Cref{lem:inner}}\geq
  \E*{\Phi_{k,T_k}} + \eta_k\E*{f_k(x_{k+1/2}) - f_k(x^*)} + \E*{\msum_{t=0}^{T_k-1}\tfrac{1}{4}\pi_k\sqn{x_{k,t+1/2} - w_{k,t}}}
  \\&\at{uses \Cref{line:xw_k+}}=
  \E*{\Phi_{k+1,0}} + \eta_k\E*{f_k(x_{k+1/2}) - f_k(x^*)} + \E*{\msum_{t=0}^{T_k-1}\tfrac{1}{4}\pi_k\sqn{x_{k,t+1/2} - w_{k,t}}}
  \\&\at{uses the previous inequality and the assumption $\E{f(\ox_k)} \in \R$}\geq
  \E*{\Phi_{k+1,0}} + \E*{\msum_{t=0}^{T_k-1}\tfrac{1}{4}\pi_k\sqn{x_{k,t+1/2} - w_{k,t}}}
  \\&
  + \eta_k\pi_k^{-1}\alpha_k^{-3}\E{f(\ox_{k+1}) - f(x^*)}
  -\eta_k\pi_k^{-1}\alpha_k^{-3}(1-\alpha_k\pi_k)\E{f(\ox_k) - f(x^*)}
\end{align*}
where \annotate.
After rearranging and using the fact that $\E{f(\ox_k)} \in \R$ and $\E{\cW_k} \in \R$, which is implied by the assumptions of the lemma, we get
\begin{align*}
  \sqn{x^*}
  &\geq
  \E{\Phi_{k,0}}
  +\eta_k\pi_k^{-1}\alpha_k^{-3}(1-\alpha_k\pi_k)\E{f(\ox_k) - f(x^*)}
  +\E{\cW_k}
  \\&\geq
  \E*{\Phi_{k+1,0}}
  + \E{\cW_{k+1}}
  + \eta_k\pi_k^{-1}\alpha_k^{-3}\E{f(\ox_{k+1}) - f(x^*)},
\end{align*}
which concludes the proof.\qed

\proofsubsection{thm:alg}

For $0 < k < n$, using \cref{eq:params}, we can obtain the following:
\begin{align*}
  \eta_k\pi_k^{-1}\alpha_k^{-3}
  -\eta_{k-1}\pi_{k-1}^{-1}\alpha_{k-1}^{-3}
  &=
  \eta (k+1)
  -\eta k
  =
  \eta_k\alpha_k^{-2}.
\end{align*}
For $k \geq n$, using \cref{eq:params}, we can obtain the following:
\begin{align*}
  \eta_k\pi_k^{-1}\alpha_k^{-3}
  -\eta_{k-1}\pi_{k-1}^{-1}\alpha_{k-1}^{-3}
  &=
  \eta n\gamma_k^{3.5}
  -\eta n\gamma_{k-1}^{3.5}
  \leq
  \tfrac{7}{2}\eta n (\gamma_k - \gamma_{k-1})\gamma_k^{2.5}
  =
  \eta\gamma_k^{2.5}
  \leq
  \eta\alpha_k\gamma_k^{3.5}
  \\&=
  \eta_k\alpha_k^{-2}.
\end{align*}
Hence, for all $k \in \N$, the following inequality holds:
\begin{align*}
  \eta_k\pi_k^{-1}\alpha_k^{-3}(1-\alpha_k\pi_k) \leq \eta_{k-1}\pi_{k-1}^{-1}\alpha_{k-1}^{-3}.
\end{align*}
Using this inequality and \Cref{lem:induction}, we can prove by induction that $T_0,\dots,T_{K-1} \in \N$ and the following inequalities hold:
\begin{align*}
  \E*{f(\ox_K) - f(x^*)}
  \leq
  \tmfrac{\alpha_{K-1}^3\pi_{K-1}}{\eta_{K-1}}\sqn{x^*}
  ,\qquad
  \E{\cW_k} \leq \sqn{x^*}.
\end{align*}
In the case where $K \geq n$ in \cref{eq:K}, we obtain the following inequality:
\begin{align*}
  \E*{f(\ox_K) - f(x^*)}
  \leq
  \tmfrac{\sqn{x^*}}{\eta n \gamma_{K-1}^{3.5}}
  \leq \epsilon.
\end{align*}
In the case where $K < n$, in \cref{eq:K}, we obtain the following:
\begin{align*}
  \sqn{x^*}/(\eta\epsilon) \leq n
  ,\quad
  K \geq n^{5/7}\brs*{\sqn{x^*}/(\eta\epsilon)}^{2/7} \geq \sqn{x^*}/(\eta\epsilon),
\end{align*}
which implies the following inequality:
\begin{align*}
  \E*{f(\ox_K) - f(x^*)}
  \leq
  \tmfrac{1}{\eta K}\sqn{x^*}
  \leq \epsilon.
\end{align*}
It remains to upper-bound $\E{\sum_{k=0}^{K-1} T_k}$ as follows:
\begin{align*}
  \E*{\msum_{k=0}^{K-1} T_k}
  &\at{is obtained similarly to the proof of \Cref{lem:inner} in \Cref{proof:lem:inner}}\leq
  K + \E*{\msum_{k=0}^{K-1} \brs*{2L\eta_k/\pi_k}^{\frac{2}{3}}\brs*{\msum_{t=0}^{T_k-1} \tfrac{1}{4}\pi_k\sqn{x_{k,t+1/2} - w_{k,t}}}^{\frac{1}{3}}}
  \\&=
  K + \E*{\msum_{k=0}^{K-1} \brs*{2L\eta_k/\pi_k}^{\frac{2}{3}}\brs*{\cW_{k+1} - \cW_k}^{\frac{1}{3}}}
  \\&\at{uses H\"older's inequality}\leq
  K + \brs*{\msum_{k=0}^{K-1} 2L\eta_k/\pi_k }^{\frac{2}{3}}\brs*{\E*{\cW_K}}^{\frac{1}{3}}
  \\&\at{uses the inequality above}\leq
  K + \brs*{\msum_{k=0}^{K-1} 2L\eta_k/\pi_k }^{\frac{2}{3}}\norm{x^*}^{\frac{2}{3}}
  \\&\at{uses \cref{eq:params}}\leq
  K + \brs*{2L\eta\norm{x^*}}^\frac{2}{3}\brs*{\msum_{k=0}^{K-1} 8\sqrt{n(k+1)}}^{\frac{2}{3}}
  \\&\leq
  K + \brs*{16\sqrt{n}L\eta\norm{x^*}}^\frac{2}{3}K
  \leq
  2K,
\end{align*}
where \annotate.\qed

\proofsubsection{cor:alg}

One can verify that the total number $C_K$ of computations of $\nabla f_i \brr{}$ and $\nabla^2 f_i\brr{}$ can be expressed as
\begin{align*}
  C_K
  =
  \cO\brr*{n + K + \msum_{k=0}^{K-1}\brs*{T_k + n\beta_k + \msum_{t=0}^{T_k-1} n \beta_{k,t}}}.
\end{align*}
Indeed, after the initial computation, we can store the full derivatives $\nabla f \brr{}$ and $\nabla^2 f\brr{}$ and update them only when $\alpha_{k}$ changes, which happens at most $\cO(1 + K/n)$ times, or when $\beta_k = 1$ or $\beta_{k,t} = 1$. These updates require $n$ computations of $\nabla f_i \brr{}$ and $\nabla^2 f_i\brr{}$.
After taking the expectation, we define $m = \min\brf{n,K}$ and obtain the following:
\begin{align*}
  \E{C_K} &= \cO\brr*{n + K + \E*{\msum_{k=0}^{K-1}\brs*{T_k + n\beta_k + \msum_{t=0}^{T_k-1} n \beta_{k,t}}}}
  \\&=
  \cO\brr*{n + K + \msum_{k=0}^{K-1}\brs*{\E{T_k} + \E*{n\beta_k} + \E*{\msum_{t=0}^{T_k-1} n \beta_{k,t}}}}
  \\&\at{uses \Cref{thm:alg}}=
  \cO\brr*{n + K + \msum_{k=0}^{K-1}\brs*{\E*{n\beta_k} + \E*{\msum_{t=0}^{T_k-1} n \beta_{k,t}}}}
  \\&\at{uses \Cref{line:beta_k}}=
  \cO\brr*{n + K + \msum_{k=0}^{K-1}\brs*{n\pi_k + \E*{\msum_{t=0}^{T_k-1} n \beta_{k,t}}}}
  \\&=
  \cO\brr*{n + K + \msum_{k=0}^{K-1}\brs*{n\pi_k + \E*{\msum_{t=0}^{+\infty} n \ind{t < T_k}\beta_{k,t}}}}
  \\&\at{use the monotone convergence theorem}=
  \cO\brr*{n + K + \msum_{k=0}^{K-1}\brs*{n\pi_k + \msum_{t=0}^{+\infty} n \E*{\ind{t < T_k}\beta_{k,t}}}}
  \\&\at{uses the fact that $\ind{t < T_k}$ is independent of $\beta_{k,t}$ and the tower rule}=
  \cO\brr*{n + K + \msum_{k=0}^{K-1}\brs*{n\pi_k + \msum_{t=0}^{+\infty} n \E*{\ind{t < T_k}\pi_k}}}
  \\&\at{use the monotone convergence theorem}=
  \cO\brr*{n + K + \msum_{k=0}^{K-1}n\pi_k \E{T_k}}
  \\&\at{use \cref{eq:params}}=
  \cO\brr*{n + K + \msum_{k=0}^{m-1}n\pi_k \E{T_k}}
  \\&\at{is implied by the proof of \Cref{lem:inner} in \Cref{proof:lem:inner}}=
  \cO\brr*{n + K + \msum_{k=0}^{m-1} n\pi_k + \E*{\msum_{k=0}^{m-1}n\pi_k\brs*{L\eta_k/\pi_k}^\frac{2}{3}
  \brs*{\msum_{t=0}^{T_k-1} \pi_k\sqn{x_{k,t+1/2} - w_{k,t}}}^{\frac{1}{3}}}}
  \\&\at{uses H\"older's inequality}\leq
  \cO\brr*{n + K + \msum_{k=0}^{m-1} n\pi_k + n(L\eta)^\frac{2}{3}\brs*{\msum_{k=0}^{m-1}\sqrt{\pi_k}
  }^{\frac{2}{3}}\brs*{\E*{\cW_m}}^{\frac{1}{3}}}
  \\&\at{is implied by the proof of \Cref{thm:alg} in \Cref{proof:thm:alg}}\leq
  \cO\brr*{n + K + \msum_{k=0}^{m-1} n\pi_k + n(L\eta\norm{x^*})^\frac{2}{3}\brs*{\msum_{k=0}^{m-1}\sqrt{\pi_k}
  }^{\frac{2}{3}}}
  \\&\at{use \cref{eq:params}}=
  \cO\brr*{K + n\log(n+1) + n(\sqrt{n}L\eta\norm{x^*})^\frac{2}{3}}
  \\&=
  \cO\brr*{K + n\log(n+1) },
\end{align*}
where \annotate.\qed

\section{Proofs for Section~\ref{sec:lb}}

\proofsubsection{lem:lb_sol}
The function $\brs{h_{d,n}(x) -\rho \<\be_{d,1},x>}$ is uniformly convex and hence has a unique solution. The exact formula for the solution given in the statement of the lemma can be verified using basic calculations.\qed

\proofsubsection{lem:lb_L}

The proof is based on the fact that each function $h_{d,n,i}(x)$ can be represented as follows:
\begin{align*}
  h_{d,n,i}(x) = \tfrac{1}{3}\norm{\mA_i x}_3^3,
\end{align*}
where $\mA_i \in \R^{\abs{\cJ_i}\times d}$ is a matrix consisting of the rows of the matrix $\mA\in\R^{d\times d}$ with indices from $\cJ_i$, where the matrix $\mA$ is defined as $\mA(t_1,\dots,t_d) = (t_1 - t_2,\dots,t_{d-1} - t_d, t_d)$.
One can verify that $\norm{\mA_i} \leq \norm{\mA} \leq 2$.  Using the fact that the Hessian $\nabla^2 \tfrac{1}{3}\norm{}_3^3$ is 2-Lipschitz, we can obtain the following upper bound:
\begin{align*}
  \tmfrac{1}{n}\msum_{i=1}^n \sqn{\nabla^2 h_{d,n,i}(x) - \nabla^2 h_{d,n,i}(x')}
  &\leq
  \tmfrac{1}{n}\msum_{i=1}^n 4 \norm{\mA_i}^4\sqn{\mA_i(x-x')}
  \\&\leq
  \tmfrac{64}{n}\msum_{i=1}^n \sqn{\mA_i(x-x')}
  \\&=
  \tmfrac{64}{n} \sqn{\mA(x-x')}
  \\&\leq
  \tmfrac{256}{n} \sqn{x-x'},
\end{align*}
which concludes the proof.\qed

\proofsubsection{thm:lb}

First, we assume that $\epsilon \leq \frac{L R^3}{300 \sqrt{n}}$ and define functions $f_i(x)$ as follows:
\begin{align*}
  f_i(x) = \ell h_{d,n,i}(x) - \ind{i = n} \cdot n\ell\rho^2\<\be_{d,1},x>,
\end{align*}
where $\ell,\rho > 0$ are defined as follows:
\begin{align*}
  \ell = \tfrac{1}{16}\sqrt{n}L
  ,\quad
  \rho = \sqrt{6/[nd(d+1)(2d+1)]}\cdot R.
\end{align*}
Let $d,K \in \N$ be defined as follows:
\begin{align*}
  d = 2K,\quad K=\floor*{\brr*{\tmfrac{L R^3}{300\sqrt{n}\epsilon}}^{2/7}}.
\end{align*}
Using \Cref{lem:lb_sol,lem:lb_L}, we can easily show that \Cref{ass:L} is satisfied and $\norm{x^*} = R$. Moreover, using \Cref{lem:lb_sol}, one can verify that the following inequality holds:
\begin{align*}
  \min_{x \in \cL_{d,K}} [f(x) - f(x^*)]
  &=
  \tfrac{2}{3}\sqrt{n}\ell \rho^3 (d-K)
  \\&=
  \tfrac{2}{3}\sqrt{n}\ell \rho^3 K
  =
  \tfrac{1}{24}nL\rho^3 K
  =
  \tmfrac{\sqrt{6} L R^3 K}{4\sqrt{n}[2K(2K+1)(4K+1)]^{3/2}}
  \\&>
  \tmfrac{L R^3}{300\sqrt{n}K^{3.5}} \geq \epsilon,
\end{align*}
where $\cL_{d,k} = \spanset{\brf{\be_{d,1},\dots,\be_{d,k}}}$ and $\cL_{d,0} = \brf{0}$.
Using standard arguments \citep[Section~4]{nesterov2021implementable}, we can show that $x_k \in \cL_{d,s_k}$, where $s_k$ is defined inductively as follows:
\begin{align*}
  s_0 = 0
  ,\quad
  s_{k+1} = s_k + \ind{i_k \equiv s_k \bmod n}.
\end{align*}
Consequently, the condition $s_k \geq K+1$ is necessary to ensure $f(x_k) - f(x^*) \leq \epsilon$. Finally, using basic properties of the Geometric distribution, we can show that
\begin{align*}
  \E{\inf\brf{k \in \N_0 : s_k = K+1}} = n(K+1),
\end{align*}
which concludes the proof in the case $\epsilon \leq \frac{L R^3}{300 \sqrt{n}}$.

In the opposite case $\epsilon > \frac{L R^3}{300 \sqrt{n}}$, we simply choose $d=1$ and define functions $f_i(x)$ as follows:
\begin{align*}
  f_i(x) =
  \begin{cases}
    \frac{2n\epsilon}{R^2}(x-R)^2 & i=1\\
    0 & i \neq 1
  \end{cases}.
\end{align*}
These functions are quadratic and obviously satisfy required properties. Moreover, since $f(0) - f(x^*) > \epsilon$, it is necessary that $i_j = 1$ for some $j < k$ to ensure $f(x_k) - f(x^*) \leq \epsilon$. Using basic properties of the Geometric distribution, we can show that $\E{\inf\brf{k \in \N : i_{k-1} = 1}} = n$, which concludes the proof.\qed

\end{document}